\documentclass[11pt,a4paper,twoside]{amsart}
\usepackage{amssymb,amsmath,amsthm}
\theoremstyle{plain}
\newtheorem{theorem}{Theorem}[section]

\newtheorem{lemma}[theorem]{Lemma}
\newtheorem{corollary}[theorem]{Corollary}
\newtheorem{proposition}[theorem]{Proposition}
\newtheorem{assumption}[theorem]{Assumption}
\theoremstyle{remark}
\newtheorem{remark}[theorem]{Remark}

\usepackage{hyperref}

\numberwithin{equation}{section}

\newcommand{\R}{\mathbb{R}}

\newcommand{\N}{\mathbb{N}}

\newcommand{\F}{\mathcal{F}}

\renewcommand{\Re}{\operatorname{Re}}
\newcommand{\I}{\infty}
\newcommand{\abs}[1]{\left\lvert #1\right\rvert}
\newcommand{\norm}[1]{\left\lVert #1\right\rVert}
\newcommand{\Lebn}[2]{\left\lVert #1 \right\rVert_{L^{#2}}}
\newcommand{\Sobn}[2]{\left\lVert #1 \right\rVert_{H^{#2}}}

\newcommand{\Jbr}[1]{\left\langle #1 \right\rangle}
\newcommand{\DE}[1]{{{#1}^\varepsilon_\delta}}
\newcommand{\DEp}[1]{{{#1}^\varepsilon_{\delta^\prime}}}
\newcommand{\IN}{\quad\text{in }}
\newcommand{\wIN}{\quad\text{weakly in }}

\def\Sch{{\mathcal S}} 
\def\({\left(}
\def\){\right)}
\def\<{\left\langle}
\def\>{\right\rangle}
\def\le{\leqslant}
\def\ge{\geqslant}

\def\Eq#1#2{\mathop{\sim}\limits_{#1\rightarrow#2}}
\def\Tend#1#2{\mathop{\longrightarrow}\limits_{#1\rightarrow#2}}

\def\d{{\partial}}
\def\g{\gamma}
\def\eps{\varepsilon}

\def\l{\lambda}

\def\si{{\sigma}}

\begin{document}
\title[Semiclassical Analysis for Hartree equation]{Semiclassical
  Analysis for Hartree equation} 
\author[R. Carles]{R\'emi Carles}
\address{CNRS \& Universit\'e Montpellier~2\\ 
    Math\'ematiques, CC~051\\
 Place Eug\`ene Bataillon\\
    34095~Montpellier~cedex, France}
\email{Remi.Carles@math.cnrs.fr}
\author[S. Masaki]{Satoshi Masaki}
\address{Department of Mathematics\\ Kyoto University\\
Kyoto 606-8502, Japan}
\email{machack@math.kyoto-u.ac.jp}
\begin{abstract}
We justify WKB analysis for Hartree equation in space dimension at
least three, in a r\'egime which is
supercritical as far as semiclassical analysis is concerned. The main
technical remark is that the nonlinear Hartree term can be considered
as a semilinear perturbation. This is in contrast with the case of the
nonlinear Schr\"odinger 
equation with a local nonlinearity, where quasilinear analysis is needed to
treat the nonlinearity. 
\end{abstract}
\maketitle
\section{Introduction}
\label{sec:intro}

We consider the semiclassical limit $\eps\to 0$ for the Hartree
equation
\begin{equation}\label{eq:r3}
  i\eps\d_t u^\eps +\frac{\eps^2}{2}\Delta u^\eps= \lambda
                (\lvert x\rvert ^{-\gamma}\ast  |u^\eps|^2)u^\eps,
                \quad \g>0,\ \l\in \R,\ x\in \R^n,
\end{equation}
in space dimension $n\ge 3$. We consider initial data of WKB type,
\begin{equation}
  \label{eq:ci}
  u^\eps(0,x) = a_0^\eps(x)e^{i\phi_0(x)/\eps},
\end{equation}
where $a_0^\eps$ typically has an asymptotic expansion as $\eps\to 0$, 
\begin{equation*}
  a_0^\eps\Eq \eps 0 a_0+\eps a_1 +\eps^2a_2+\ldots,\quad a_j\text{
  independent of }\eps\in ]0,1].
\end{equation*}
The approach that we follow is closely related to the pioneering works
of P.~G\'erard \cite{PGX93}, and E.~Grenier \cite{Grenier98}, for the
nonlinear Schr\"odinger equation with local nonlinearity:
\begin{equation}
  \label{eq:nls}
  i\eps\d_t u^\eps +\frac{\eps^2}{2}\Delta u^\eps = f\(\lvert
  u^\eps\rvert^2\)u^\eps, 
\end{equation}
where the function $f$ is smooth, and real-valued. In
\cite{Grenier98}, the assumption $f'>0$ is necessary for the arguments
of the proof. More recently, this assumption was
relaxed in \cite{AC-BKW}, allowing to consider the case $f(y)=+y^\si$, $\si \in
\N$. Moreover, it is noticed in \cite{AC-BKW} that to carry out a WKB
analysis in Sobolev spaces for \eqref{eq:nls}, the assumption $f'\ge
0$ is essentially necessary. Typically, in the case $f'<0$, working with
analytic data is necessary, and sufficient as shown in
\cite{PGX93,ThomannAnalytic}. The reason is that the local
nonlinearity is analyzed through \emph{quasilinear} arguments in
\cite{Grenier98,AC-BKW}, and $f'$ determines the velocity
of a wave equation: if $f'>0$, then the wave equation is hyperbolic, and
$f'<0$, the underlying operator becomes elliptic. 
\smallbreak

The above discussion is altered for the Hartree type
nonlinearity. Typically, no assumption is made on the sign of $\l$
here. As noticed in \cite{AC-SP} in the special case of the
Schr\"odinger--Poisson system, the nonlocal nonlinearity in
\eqref{eq:r3} can be handled by \emph{semilinear} arguments. However,
a quasilinear analysis is needed to handle the convective coupling.
\smallbreak

There are at least to motivations to study this question, besides the
general picture of justifying  approximations motivated by physics. As
remarked in \cite{CaJHDE} in the case of a local nonlinearity, WKB
analysis and a geometrical transform 
can help understand the behavior of a wave function near a focal
point, in a supercritical r\'egime. In \cite{SatoshiAHP}, other
informations were obtained thanks to a different approach, in the case
of a Hartree type nonlinearity. The approach of \cite{SatoshiAHP} and
the results of the present paper will certainly be helpful to improve
the understanding of the focusing phenomenon in semiclassical
analysis. Another application of the WKB analysis for \eqref{eq:r3}
concerns the Cauchy problem for the Hartree equation, that is
\eqref{eq:r3} with $\eps=1$. Following the approach initiated in
\cite{BGTMRL,BGTENS,CCT,CCT2,Lebeau01,Lebeau05}, we can prove an
ill-posedness result, together with a loss of regularity; see
Corollary~\ref{cor:loss}. 
\begin{assumption}\label{asmp:existence}
Let $n \ge 3$ and $\max(n/2-2,0) < \gamma \le n-2$.
We suppose the following conditions with some $s>n/2+1$:\\
$\bullet$ The initial amplitude $a_0^\varepsilon \in H^s(\R^n)$, uniformly for
  $\varepsilon \in [0,1]$: 
there exists a constant $C$ independent of $\varepsilon$
  such that $\Sobn{a_0^\varepsilon}{s} \le C$. \\
$\bullet$  The initial phase $\phi_0$ satisfies
\begin{equation*}
  \lvert \phi_0(x)\rvert + \lvert \nabla \phi_0(x)\rvert \Tend
  {|x|}\infty 0, 
\end{equation*}
and $\nabla^2\phi_0 \in H^s(\R^n)$.
Moreover, there exists $q_0 \in ]n/(\gamma+1),n[$ such that $\nabla
\phi_0 \in L^{q_0}(\R^n)$. 
\end{assumption}
\begin{remark}
  We will see that the above assumption implies that $\phi_0$ and
  $\nabla \phi_0$ are bounded, and enjoy some extra integrability
  properties. See Remark~\ref{rem:hyp}. 
\end{remark}
\begin{remark}
  By employing the geometrical reduction made in \cite{AC-SP},
  we can relax the assumption $\lvert \phi_0(x)\rvert + \lvert \nabla 
  \phi_0(x)\rvert \to 0$ as $|x|\to \infty$ in a sense. 
  Indeed, we can replace the initial phase $\phi_0$ with $\phi_0 + 
  \phi_{\text{quad}}$, where $\phi_{\text{quad}} \in C^\I(\R^n)$
  is a polynomial of degree at most two.
\end{remark}
For $s>n/2$, we denote by $X^s(\R^n)$ the Zhidkov space 
\[
        X^s(\R^n) = \{ u \in L^\I(\R^n)| \nabla u \in H^{s-1}(\R^n) \}.
\]
This space was introduced in \cite{Zhidkov} (see also
\cite{ZhidkovLNM}) in the case 
$n=1$, and its study was generalized to the multidimensional case in
\cite{Gallo}. We denote
\[
        \norm{u}_{X^s} := \Lebn{u}{\I} + \Sobn{\nabla u}{s-1}.
\]
We write $H^s = H^s(\R^n)$ and $X^s = X^s(\R^n)$.
\begin{theorem}\label{thm:existence}
Let Assumption \ref{asmp:existence} be satisfied.
There exists $T>0$ independent of $\varepsilon$ and $s>1+n/2$, 
and a unique solution $u^\varepsilon \in C([0,T];H^s)$ to the equation
\eqref{eq:r3}--\eqref{eq:ci}. Moreover, it can be written in 
 the form $u^\eps = a^\eps e^{i\phi^\eps/\eps}$, where
$a^\eps$ is complex-valued, $\phi^\eps$ is real-valued,  with
\begin{align*}
  a^\varepsilon \in C([0,T]; H^s)&,\ \phi^\varepsilon \in C\([0,T];
  L^\I \cap L^{\frac{nq_0}{n-q_0}}\),\\ 
\text{ and }&\nabla
  \phi^\varepsilon \in C\([0,T]; X^{s+1} \cap L^{q_0}\). 
\end{align*}
Moreover, if $\phi_0 \in L^{p_0}$ for some $p_0\in]n/\gamma,nq_0/(n-q_0)[$
then $\phi^\varepsilon \in C([0,T]; L^{p_0})$.
\end{theorem}
Note that obtaining a local existence time $T$ which is independent of
$\eps\in ]0,1]$ is already a non-trivial information, at least for a
focusing nonlinearity $\l<0$. Using classical
results on the Cauchy problem for Hartree equation \cite{GV82}, and a
scaling argument, would yield an existence time that goes to zero with
$\eps$. 
Taking $q_0=2$, we immediately obtain the following corollary:
\begin{corollary}
Let Assumption~\ref{asmp:existence} be satisfied.
Let $u^\varepsilon = a^\varepsilon
e^{i\phi^\varepsilon/\varepsilon}$ be the solution given in
Theorem~\ref{thm:existence}. 
If $\gamma \in ]n/2-1,n-2]$ and $\nabla \phi_0 \in H^{s+1}$,
then 
\begin{equation*}
  \phi^\varepsilon \in C\([0,T];X^{s+2}\cap
L^{\frac{2n}{n-2}}\).
\end{equation*}
\end{corollary}
With this local existence result, we can justify a WKB expansion,
provided that the initial data have a suitable expansion as $\eps\to
0$. 
\begin{assumption}\label{asmp:expansion}
Let $N$ be a positive integer. 
We suppose Assumption~\ref{asmp:existence} with some $s>n/2 + 2N +1$.
Moreover, the initial amplitude $a_0^\varepsilon$ writes
\begin{equation}
        a_0^\varepsilon = a_0 + \sum_{j=1}^N \varepsilon^j a_j +
        \varepsilon^N r^\varepsilon_N, 
\end{equation}
where $a_j \in H^{s}$ ($0 \le j \le N$) and
$\Sobn{r^\varepsilon_N}{s}\to 0$ as $\varepsilon\to 0$.  
\end{assumption}

\begin{theorem}\label{thm:expansion}
Let Assumption \ref{asmp:expansion} be satisfied.
Let $u^\varepsilon=a^\varepsilon e^{i\phi^\varepsilon/\varepsilon}$ be
the unique solution given in Theorem~\ref{thm:existence}. 
Then, there exist $(b_j,\varphi_j)_{0\le j\le N}$, with
\begin{align*}
  &b_j \in
C([0,T];H^{s-2j}),\\
&\varphi_j\in
C\([0,T];L^\I \cap L^{\frac{nq_0}{n-q_0}}\)\text{ and }\nabla \varphi_j \in
C\([0,T];X^{s-2j+1}\cap L^{q_0}\), 
\end{align*}
such that:
\begin{align*}
        a^\varepsilon &= b_0 + \sum_{j=1}^{N} \varepsilon^j b_j
        +o\(\varepsilon^N\) \IN C([0,T];H^{s-2N}), \\ 
        \phi^\varepsilon &= \varphi_0 + \sum_{j=1}^{N-1} \varepsilon^j
        \varphi_j +o\(\varepsilon^N\) \IN C([0,T];L^\I \cap
        L^{\frac{nq_0}{n-q_0}}), \\ 
        \nabla\phi^\varepsilon &= \nabla\varphi_0 + \sum_{j=1}^{N}
        \varepsilon^j \nabla\varphi_j +o\(\varepsilon^N\) \IN
        C\([0,T];X^{s-2N+1} \cap L^{q_0}\).
\end{align*}
Moreover, for $j \ge 1$, $\varphi_j \in L^p$ for all $p>n/\gamma$,
and $\nabla \varphi_j 
\in L^q$ for all $q>n/(\gamma+1)$. 
\end{theorem}
\begin{corollary}\label{cor:bkw}
 Let Assumption \ref{asmp:expansion} be satisfied.
The solution $u^\varepsilon$ given in Theorem~\ref{thm:existence} has
the following asymptotic expansion, as $\eps \to 0$:
\begin{equation*}
  u^\eps = e^{i\varphi_0/\eps}\( \beta_0+\eps \beta_1+\ldots
  +\eps^{N-1}\beta_{N-1}+\eps^{N-1}\rho^\eps\), \quad \lVert
  \rho^\eps\rVert_{H^{s-2N+2}} \Tend \eps 0 0 ,
\end{equation*}
where $\varphi_0$ is given by Theorem~\ref{thm:expansion}, and $\beta_j
\in C([0,T];H^{s-2j})$ is a smooth function of
$(b_k,\varphi_{k+1})_{0\le k\le j}$. For instance,
\begin{equation*}
  \beta_0 = b_0 e^{i\varphi_1}\quad ;\quad \beta_1 = b_1
  e^{i\varphi_1} +i\varphi_2 b_0 e^{i\varphi_1}. 
\end{equation*}
\end{corollary}
\begin{corollary}\label{cor:loss}
Let $n\ge 5$, $\l\in \R$, $\max(n/2-2,2)<\g \le n-2$ and $0<s<s_c=
\g/2-1$. There exists a sequence of 
initial data
\begin{equation*}
  (\psi_0^h)_{0<h\le 1},\ \psi_0^h\in \Sch(\R^n),\quad \lVert
  \psi_0^h \rVert_{H^s}\Tend h 0 0,
\end{equation*}
a sequence of times $t^h\to 0$, such that the solution to
\begin{equation*}
  i\d_t \psi^h +\frac{1}{2}\Delta \psi^h = \l \(\lvert
  x\rvert^{-\g}\ast \lvert \psi^h\rvert^2\)\psi^h \quad ;\quad
  \psi^h_{\mid t=0}= \psi_0^h
\end{equation*}
satisfies
\begin{equation*}
  \left\lVert \psi^h(t^h)\right\rVert_{H^k}\Tend h 0 +\infty, \quad
  \forall k>\frac{s}{1+s_c-s}=\frac{s}{\g/2-s}. 
\end{equation*}
\end{corollary}
Note that unlike in the case of Schr\"odinger equations with local
nonlinearity, considering a large space dimension is necessary to
observe this phenomenon: in low space dimensions, Hartree equations
are locally well-posed in Sobolev spaces of positive regularity (see
e.g. \cite{CazCourant,GV82}). 

Using Sobolev embedding, one could infer
a loss of regularity at the level of the energy space (consider $s>1$
and $\frac{s}{\g/2-s}=1$, hence $\g=4s>1$), in the spirit of
\cite{Lebeau05} (see also \cite{AC-perte,ThomannAnalytic} for
Schr\"odinger equations), provided that the space dimension is $n\ge
7$. 
\smallbreak

The rest of this paper is organized as follows. In the next paragraph,
we present the general strategy adopted in this paper. In \S\ref{sec:prelim},
we collect some technical estimates. Theorem~\ref{thm:existence} is
proved in \S\ref{sec:existence}, and Theorem~\ref{thm:expansion} is
proved in  \S\ref{sec:expansion}, as well as
Corollary~\ref{cor:bkw}. Finally, Corollary~\ref{cor:loss} is inferred
in \S\ref{sec:loss}.

\section{General strategy}
\label{sec:strategy}

To prove Theorem~\ref{thm:existence} and Theorem~\ref{thm:expansion},
we follow the same strategy as in \cite{Grenier98}. 
Seek a solution $u^\varepsilon$ to \eqref{eq:r3}--\eqref{eq:ci}
represented as 
\begin{equation}\label{eq:u-ap}
        u^\varepsilon(t,x)=a^\varepsilon(t,x)
        e^{i\phi^\varepsilon(t,x)/\varepsilon},
\end{equation}
with a complex-valued space-time function $a^\varepsilon$ and a
real-valued space-time function $\phi^\varepsilon$. Note that $a^\eps$
is expected to be complex-valued, even if its initial value $a_0^\eps$
is real-valued. 
We remark that the phase function $\phi^\varepsilon$ also depends on
the parameter $\eps$. 
Substituting the form \eqref{eq:u-ap} into \eqref{eq:r3}, we obtain
\begin{align*}
        i \varepsilon \left(\partial_t a^\varepsilon + a^\varepsilon i
        \frac{\partial_t \phi^\varepsilon}{\varepsilon}\right) 
        + \frac{\varepsilon^2}{2} \left( \Delta a^\varepsilon +
        2\left(\nabla a^\varepsilon \cdot i  \frac{\nabla
        \phi^\varepsilon}{\varepsilon}\right)  
        - a^\varepsilon \frac{|\nabla
        \phi^\varepsilon|^2}{\varepsilon^2} + a^\varepsilon i
        \frac{\Delta \phi^\varepsilon}{\varepsilon}\right) 
        \\= \lambda (|x|^{-\gamma}*|a^\varepsilon|^2)a^\varepsilon.
\end{align*}
To obtain a solution of the above equation (hence, of \eqref{eq:r3}),
we choose to consider the following system:
\begin{align}\label{eq:originalsystem1}
        \begin{cases}
                \partial_t a^\varepsilon + \nabla a^\varepsilon \cdot
                \nabla \phi^\varepsilon + \dfrac{1}{2}a^\varepsilon
                \Delta \phi^\varepsilon 
                = i \dfrac{\varepsilon}{2}\Delta a^\varepsilon \\
                \partial_t \phi^\varepsilon + \dfrac{1}{2}|\nabla
                \phi^\varepsilon|^2 + \lambda
                (|x|^{-\gamma}*|a^\varepsilon|^2) =0. 
        \end{cases}
\end{align}
This choice is essentially the same as the one introduced by
E.~Grenier \cite{Grenier98}. 
We consider this with the initial data
\begin{align}\label{eq:originalsystem2}
        a^\varepsilon_{|t=0} &=a_0^\varepsilon, &
        \phi^\varepsilon_{|t=0} &=\phi_0. 
\end{align}
From now on, we work only on
\eqref{eq:originalsystem1}--\eqref{eq:originalsystem2}. We first prove
that it admits a unique solution with suitable regularity (see
Theorem~\ref{thm:existence}), hence providing a solution to
\eqref{eq:r3}--\eqref{eq:ci}. The asymptotic expansion
(Theorem~\ref{thm:expansion}) then follows by the similar arguments. 
\smallbreak

To conclude this paragraph, we remark that uniqueness for
\eqref{eq:r3}--\eqref{eq:ci} in the class $C([0,T];H^s)$, $s>n/2+1$,
is a straightforward consequence of \cite{GV82} (see also
\cite{CazCourant}) in space dimension
$3\le n\le 5$, since the parameter
$\eps\in ]0,1]$ can be considered fixed. Indeed, in that case, one has
$x\mapsto |x|^{-\g}\in L^p+L^\infty$ for some $p\ge 1$ and $p>n/4$,
since $\g \le n-2<4$. Uniqueness is actually obtained in the
weaker class of finite energy solutions. Since we work at a higher
degree of regularity, we can simply notice that the nonlinear potential
$|x|^{-\g}\ast |u^\eps|^2$ is bounded in $L^\infty([0,T]\times\R^n)$:
for $\chi \in C_0^\infty\(\R^n;[0,1]\)$, $\chi=1$ near $x=0$, we have
from Young's inequality,
\begin{align*}
  \left\lVert |x|^{-\g}\ast |u^\eps|^2\right\rVert_{L^\infty(\R^n)}
  &\le \left\lVert \chi |x|^{-\g}\right\rVert_{L^1(\R^n)} \left\lVert
  u^\eps\right\rVert_{L^\infty(\R^n)}^2\\
& + \left\lVert (1-\chi)
  |x|^{-\g}\right\rVert_{L^\infty(\R^n)} \left\lVert 
  u^\eps\right\rVert_{L^2(\R^n)}^2.
\end{align*}
Since we work with an $H^s$ regularity, $s>n/2+1$, the above right
hand side is bounded, and uniqueness follows
from standard energy estimates in $L^2$.

\section{Preliminary estimates}\label{sec:prelim}
We first recall a consequence of the Hardy-Littlewood-Sobolev inequality,
which can be found in \cite[Th.~4.5.9]{Hormander1} or
\cite[Lemma~7]{PG05}: 
\begin{lemma}\label{lem:Zhidkov1}
If $\varphi \in \mathcal{D}^\prime(\R^n)$ is such that $\nabla \varphi
\in L^p (\R^n)$ for some $p\in]0,n[$, 
then there exists a constant $\gamma$ such that $\varphi-\gamma \in
L^q(\R^n)$, with $1/p=1/q+1/n$. 
\end{lemma}
\begin{remark}
  The limiting case $\g = n-2$ corresponds to the
  Schr\"odinger--Poisson system considered in \cite{AC-SP}, with
  suitable conditions at infinity to integrate the Poisson equation. 
\end{remark}
\begin{remark}\label{rem:hyp}
By Lemma \ref{lem:Zhidkov1} and Sobolev inequality,
Assumption~\ref{asmp:existence} implies  $\phi_0 \in
L^{\frac{nq_0}{n-q_0}}\cap L^\I$,
and $\nabla \phi_0 \in L^{q_0} \cap X^{s+1}$.
Note that $2n/(n-2)<n$ if $n \ge 5$.
Therefore, in this case, we can always find $q_0$ in
$]n/(\gamma+1),n[$ such that $\nabla \phi_0 \in L^{q_0}$. 
\end{remark}
The next two lemmas can be found in \cite{lannesJFA}: 
\begin{lemma}[Commutator estimate]
Let $s \ge 0$ and $1<p<\I$. Set $\Lambda = (1-\Delta)^{1/2}$.
Then, it holds  that 
\[
        \Lebn{\Lambda^s(fg) - f \Lambda^s g}{p} \le
        c(\Lebn{\nabla f}{\I} \Lebn{\Lambda^{s-1}g}{p} +
        \Lebn{\Lambda^s f}{p}\Lebn{g}{\I}). 
\]
\end{lemma}
%
%
%
%
\begin{lemma}
Let $s >0$ and $1<p<\I$. There exists $C>0$ such that
\[
        \Lebn{\Lambda^s(fg)}{p} \le C(\Lebn{\Lambda^s f}{p}
        \Lebn{g}{\I} + \Lebn{f}{\I}\Lebn{\Lambda^s g}{p}),\quad
        \forall f,g \in W^{s,p} \cap L^\I.
\]
\end{lemma}
The following lemma is crucial for our analysis:
\begin{lemma}\label{lem:Hartree}
Let $n \ge 3$, $k\ge 0$, and $s_1,s_2 \in \R$.
Let $\gamma>0$ satisfying $n/2-k < \gamma \le n-k-s_1+s_2$.
Then, there exists $C_s$ such that 
\[
        \norm{\lvert\nabla\rvert^k(|x|^{-\gamma} * f)}_{H^{s_1}} \le C_s
        (\norm{f}_{H^{s_2}} + \norm{f}_{L^1}),\quad \forall f \in L^1
        \cap H^{s_2}. 
\]
\end{lemma}
\begin{proof}
Since $\F |x|^{-\gamma} = C|\xi|^{-n+\gamma}$, it holds that
\[
        \norm{|\nabla|^k(|x|^{\gamma}*f)}_{H^{s_1}} 
        = C\norm{\Jbr{\xi}^{s_1} |\xi|^{-n+\gamma+k} \F f}_{L^2}.
\]
The high frequency part $(|\xi|>1)$ is bounded by $C\norm{f}_{H^{s_2}}$
if $-n+\gamma+k +s_1-s_2\le 0$.
On the other hand, the low frequency part $(|\xi| \le 1)$ is bounded by
\begin{align*}
        C\Lebn{\F f}{\I}\int_{|\xi| \le 1}|\xi|^{2(-n+\gamma+k)}d\xi
        \le C\Lebn{f}{1} 
\end{align*}
if $2(-n+\gamma+k)>-n$, that is, if $\gamma>n/2-k$.
\end{proof}

\section{Existence result: proof of Theorem~\ref{thm:existence}}
\label{sec:existence}
Operating $\nabla$ to the equation for $\phi^\eps$ in
\eqref{eq:originalsystem1} and putting $v^\varepsilon:=\nabla
\phi^\varepsilon$, we obtain the following system:
\begin{equation}\label{eq:system1}
\left\{
  \begin{aligned}
    &\partial_t a^\varepsilon + v^\varepsilon \cdot \nabla
        a^\varepsilon + \dfrac{1}{2}a^\varepsilon \nabla\cdot
        v^\varepsilon 
        = i \dfrac{\varepsilon}{2} \Delta a^\varepsilon,&&
        a^\varepsilon_{|t=0} = a_0^\varepsilon, \\ 
        &\partial_t v^\varepsilon + v^\varepsilon \cdot \nabla
        v^\varepsilon + \lambda \nabla \(|x|^{-\gamma}
        \ast |a^\varepsilon|^2\) = 0, && v^\varepsilon_{|t=0}=\nabla \phi_0.
  \end{aligned}
\right.
\end{equation}

We first construct the solution $(a^\varepsilon,v^\varepsilon)$ to the
system \eqref{eq:system1}. 
\begin{proposition}\label{prop:existence-sys}
Let Assumption~\ref{asmp:existence} be satisfied.
There exists $T>0$ independent of $\varepsilon$ and $s$,
such that for all $\eps\in [0,1]$, \eqref{eq:system1} has a
unique solution 
\begin{equation*}
  (a^\varepsilon,v^\varepsilon)\in C\([0,T]; H^s\times \(X^{s+1}\cap
L^{\frac{nq_0}{n-q_0}}\)\).
\end{equation*}
Moreover, the norm of $ (a^\varepsilon,v^\varepsilon)$ is bounded
uniformly for $\eps \in ]0,1]$. 
\end{proposition}
\subsection{Regularized system}
We shall prove the existence of the solution to the system \eqref{eq:system1}
by taking the limit of the solutions to the corresponding regularized system.
We take $\varphi \in C_0^\I(\R^n)$, with
$\int_{\R^n} \varphi(x) dx = 1$ and $\varphi \ge 0$ and set 
\begin{equation}\label{def:Jd}
        J_\delta f = \varphi_\delta * f
\end{equation}
where $\varphi_\delta = \delta^{-n}\varphi(x/\delta)$.
We first treat the following regularized system:
\begin{equation}\label{eq:regsys1}
\left\{
        \begin{aligned}
        &\partial_t a^\varepsilon_\delta + J_\delta (
        v^\varepsilon_\delta \cdot \nabla J_\delta
        a^\varepsilon_\delta) + 
        \dfrac{1}{2} a^\varepsilon_\delta \nabla\cdot J_\delta
        v^\varepsilon_\delta = i \dfrac{\varepsilon}{2} \Delta
        J_\delta^2 a^\varepsilon_\delta\ ;&& a^\varepsilon_{\delta|t=0}
        = a_0^\varepsilon.\\  
        &\partial_t v^\varepsilon_\delta + J_\delta (
        v^\varepsilon_\delta \cdot \nabla J_\delta
        v^\varepsilon_\delta) + \lambda \nabla J_\delta (|x|^{-\gamma}
        *| a^\varepsilon_\delta|^2) = 0\ ;&&
        v^\varepsilon_{\delta|t=0}=\nabla \phi_0. 
        \end{aligned}
\right.
\end{equation}
The point is that the regularized equations \eqref{eq:regsys1}  have
been chosen so that 
the Cauchy problem can be solved as in the standard framework of
Sobolev and Zhidkov spaces: 
\begin{lemma}\label{lem:regexistence}
Let Assumption \ref{asmp:existence} be satisfied.
For all $\varepsilon \in [0,1]$ and $\delta\in]0,1]$, there exists
$T_\delta^\varepsilon>0$ such that 
the Cauchy problem \eqref{eq:regsys1} has a unique solution
$(a^\varepsilon_\delta,v^\varepsilon_\delta)\in
C^1([0,T_\delta^\varepsilon],H^{s+1}\times X^{s+2}\cap
L^{\frac{2n}{n-2}})$. 
\end{lemma}
\begin{proof}
The proof is based on the usual theorem for ordinary differential equations.
We use the following estimates
\begin{align*}
        \norm{J_\delta( v^\varepsilon_\delta\cdot \nabla J_\delta
                a_\delta^\varepsilon)}_{H^{s+1}} 
                &\le
                C\norm{v^\varepsilon_\delta}_{H^{s+1}}\norm{\nabla
                J_\delta a_\delta^\varepsilon}_{H^{s+1}} \\ 
                &\le C \delta^{-1}
                \norm{v^\varepsilon_\delta}_{H^{s+2}}
                \norm{a^\varepsilon_\delta}_{H^{s+1}},    
\end{align*}
and
\begin{align*}
        \norm{ a^\varepsilon_\delta \nabla \cdot J_\delta
                v_\delta^\varepsilon}_{H^{s+1}} 
                &\le C\delta^{-1}
                \norm{a^\varepsilon_\delta}_{H^{s+1}}
                \norm{v^\varepsilon_\delta}_{X^{s+2}},\\ 
        \norm{\Delta J_\delta^2 a_\delta^\varepsilon}_{H^{s+1}}
                &\le C \delta^{-2} \norm{a_\delta^\varepsilon}_{H^{s+1}}, \\
        \norm{J_\delta( v^\varepsilon_\delta\cdot \nabla J_\delta
                v_\delta^\varepsilon)}_{X^{s+2}} 
                &\le C \delta^{-1}\norm{v^\varepsilon_\delta}_{X^{s+2}}^2, \\
        \norm{\nabla J_\delta (|x|^{-\gamma}*|
                a_\delta^\varepsilon|^2)}_{X^{s+2}} 
                &\le C \norm{\Delta (|x|^{-\gamma}*|
                a_\delta^\varepsilon|^2)}_{H^{s+1}}\\ 
                &\le C \norm{a_\delta^\varepsilon}_{H^{s+1}}^2.
\end{align*}
We have applied Lemma \ref{lem:Hartree} with $k=2$ and $s_1=s_2=s+1$.
We note that the space $X^{s+2}\cap L^{\frac{2n}{n-2}}$ with norm
$\norm{\cdot}_{X^{s+2}}$ is complete. 
\end{proof}
\subsection{Uniform bound}
We shall establish an upper bound for the $H^s$ norm and $X^{s+1}$ norm of
$a_\delta^\varepsilon$ and $v^\varepsilon_\delta$ for $s>n/2+1$,
respectively. 
We first estimate the $H^s$ norm of $a^\varepsilon_\delta$.
We use the following convention for the scalar product in $L^2$:
\[
        \Jbr{\varphi, \psi} := \int_{\R^n} \varphi(x) \overline{\psi(x)} dx.
\]
Set $\Lambda=(I-\Delta)^{1/2}$. We shall estimate
\[
        \frac{d}{dt} \Sobn{\DE{a}}{s}^2 = 2\Re \Jbr{\partial_t
        \Lambda^s \DE{a}, \Lambda^s \DE{a}}. 
\]
Since $[\Lambda^s,\nabla]=0$ and $[\Lambda^s,J_\delta]=0$, by
commuting $\Lambda^s$ with the equation for $\DE{a}$, we find: 
\begin{equation}\label{eq:unifbd1}
        \partial_t \Lambda^s \DE{a} + J_\delta \Lambda^s  ( \DE{v}
        \cdot \nabla J_\delta \DE{a}) 
        + \frac{1}{2} \Lambda^s  ( \DE{a} \nabla J_\delta \cdot \DE{v})
         -i \frac{\varepsilon}{2} J_\delta \Delta  J_\delta \Lambda^s
        \DE{a}= 0. 
\end{equation}
The coupling of the second term and $\Lambda^s \DE{a}$ is written as
\begin{align*}
        \Jbr{\Lambda^s( \DE{v}\cdot \nabla J_\delta \DE{a}),  J_\delta
        \Lambda^s \DE{a} } =&  
        \Jbr{ \DE{v}\cdot \nabla J_\delta \Lambda^s \DE{a}, J_\delta
        \Lambda^s \DE{a} } \\ 
        &+ \Jbr{[\Lambda^s, \DE{v}]\cdot  \nabla J_\delta \DE{a},
        J_\delta \Lambda^s  \DE{a} }, 
\end{align*} 
where we have use the fact that $\Jbr{J_\delta f,g} = \Jbr{f,J_\delta
  g}$ for any $f$ and $g$. We see from the integration by parts that
\begin{equation}\label{eq:unifbd2}
        |\Re \Jbr{ \DE{v}\cdot \nabla J_\delta \Lambda^s  \DE{a},
        J_\delta  \Lambda^s  \DE{a} }| 
        \le \frac{1}{2} \Lebn{\nabla \DE{v}}{\I} \Lebn{ J_\delta
        \Lambda^s \DE{a}}{2}^2. 
\end{equation}
Moreover, the commutator estimate shows that
\begin{multline}\label{eq:unifbd3}
        |\Re \Jbr{[\Lambda^s, \DE{v}]\cdot  \nabla J_\delta \DE{a},
        J_\delta \Lambda^s  \DE{a} }|\\ 
        \le C(\Sobn{ \nabla \DE{v}}{s-1} \Lebn{J_\delta \nabla
        \DE{a}}{\I} + \Lebn{ \nabla \DE{v}}{\I} \Sobn{J_\delta \nabla
        \DE{a}}{s-1} ) \Lebn{J_\delta \Lambda^s \DE{a}}{2} 
\end{multline}
We estimate the third term of \eqref{eq:unifbd1} by the Kato--Ponce
inequality as 
\begin{multline}\label{eq:unifbd4}
        |\Re \Jbr{ \Lambda^s  ( \DE{a} \nabla J_\delta \cdot \DE{v}),
        \Lambda^s \DE{a}}|\\ 
        \le C( \Lebn{ \DE{a}}{\I} \Sobn{ J_\delta \nabla \DE{v}}{s}  +
        \Sobn{ \DE{a}}{s} \Lebn{J_\delta \nabla \DE{v}}{\I}  ) \Lebn{
        \Lambda^s \DE{a}}{2} 
\end{multline}
and the last term vanishes since
\begin{equation}\label{eq:unifbd5}
        \Re\Jbr{-i \Delta  J_\delta \Lambda^s \DE{a}, J_\delta
        \Lambda^s \DE{a}} = \Re i \Lebn{\nabla J_\delta \Lambda^s
        \DE{a}}{2}^2 = 0. 
\end{equation} 
Therefore, summarizing \eqref{eq:unifbd1}--\eqref{eq:unifbd5}, we end up with 
\begin{equation*}
        \frac{d}{dt} \Sobn{\DE{a}}{s}^2 \le C (\norm{\DE{a}
        }_{W^{1,\infty}} + \Lebn{\nabla \DE{v}}{\I})( \Sobn{\DE{a}}{s}
        + \norm{\DE{v}}_{X^{s+1}})\Sobn{\DE{a}}{s}, 
\end{equation*}
hence
\begin{equation}\label{eq:unifbd6}
        \frac{d}{dt} \Sobn{\DE{a}}{s}^2 \le C (\norm{\DE{a}
        }_{W^{1,\infty}} + \Lebn{\nabla \DE{v}}{\I})( \Sobn{\DE{a}}{s}^2
        + \norm{\DE{v}}_{X^{s+1}}^2). 
\end{equation}
Let us proceed to the estimate of $\DE{v}$. We denote the operator
$\Lambda^s \nabla$ by $Q$. 
From the equation for $\DE{v}$, we have
\begin{equation}\label{eq:unifbd7}
        \partial_t Q \DE{v} + J_\delta Q ( \DE{v} \cdot \nabla
        J_\delta \DE{v}) + Q \nabla \Lambda^{s}J_\delta
        (|x|^{-\gamma}*| \DE{a}|^2)=0 
\end{equation}
We consider the coupling of this equation and $Q\DE{v}$.
The second term can be written as
\begin{align*}
        \Jbr{Q ( \DE{v} \cdot \nabla J_\delta \DE{v}), J_\delta Q \DE{v}}
        =& \Jbr{ \DE{v}\cdot \nabla J_\delta Q \DE{v}, J_\delta Q \DE{v}} \\
        &+ \Jbr{[Q, \DE{v}]\cdot \nabla J_\delta \DE{v}, J_\delta Q
        \nabla\DE{v}}. 
\end{align*}
As the previous case,  integration by parts shows
\begin{equation}\label{eq:unifbf8}
        |\Re \Jbr{ \DE{v}\cdot \nabla J_\delta Q\DE{v}, J_\delta Q \DE{v}}|
        \le \frac{1}{2} \Lebn{ \nabla \cdot \DE{v}}{\I} \Lebn{J_\delta
        Q \DE{v}}{2}, 
\end{equation}
and the commutator estimate also shows
\begin{multline}\label{eq:unifbd9}
        |\Re \Jbr{[Q, \DE{v}]\cdot \nabla J_\delta \DE{v}, J_\delta Q
        \DE{v}}|\\ 
        \le C (\Sobn{ \nabla\DE{v}}{s} \Lebn{J_\delta \nabla
        \DE{v}}{\I} + \Lebn{ \nabla \DE{v}}{\I}\Sobn{J_\delta \nabla
        \DE{v}}{s}) 
        \Lebn{J_\delta Q \DE{v}}{2}.
\end{multline}
For the estimate of the Hartree nonlinearity, we use Lemma
\ref{lem:Hartree} with $k=2$ and $s_1=s_2=s$, to obtain 
\begin{align}
        \Lebn{\lambda J_\delta \Lambda^s  \nabla^2 (|x|^{-\gamma}*|
        \DE{a}|^2)}{2} 
        &\le C (\Sobn{| \DE{a}|^2}{s} + \Lebn{| \DE{a}|^2}{1}) \nonumber\\
        &\le C (\Sobn{ {\DE{a}} }{s} \Lebn{ {\DE{a}} }{\I} + \Lebn{
        {\DE{a} }}{2}^2). \label{eq:unifbd10} 
\end{align}
Summarizing \eqref{eq:unifbd7}--\eqref{eq:unifbd10}, we deduce that
\begin{equation*}
        \frac{d}{dt} \Sobn{\nabla\DE{v}}{s}^2 \le C (\Lebn{\DE{a} }{2}
        + \Lebn{\DE{a} }{\I} + \Lebn{\nabla \DE{v}}{\I})(
        \Sobn{\DE{a}}{s} +
        \Sobn{\nabla\DE{v}}{s})\Sobn{\nabla\DE{v}}{s}, 
\end{equation*}
and so that 
\begin{equation}\label{eq:unifbd11}
        \frac{d}{dt} \Sobn{\nabla\DE{v}}{s}^2 \le C (\Lebn{\DE{a} }{2} +
        \Lebn{\DE{a} }{\I} + \Lebn{\nabla \DE{v}}{\I})(
        \Sobn{\DE{a}}{s}^2 +\Sobn{\nabla\DE{v}}{s}^2). 
\end{equation}
Using Lemma~\ref{lem:Zhidkov1}, we see that the above estimate yields
an $L^{2n/(n-2)}$ estimate for $\DE{v}$. Interpolating with a suitable
$\dot H^k$ norm shows that the $L^\infty$ norm of  $\DE{v}$ is
estimated as above. Alternatively, 
integrating the second equation of \eqref{eq:regsys1} with respect to
time, Sobolev embedding directly yields a similar estimate for $\DE{v}$ in
$L^\infty(\R^n)$. 

Now, putting $\DE{M}(t):= \Sobn{\DE{a}(t)}{s}^2 +
\norm{\nabla\DE{v}(t)}_{H^{s}}^2+\norm{\DE{v}(t)}_{L^\infty}^2 $, we
conclude from \eqref{eq:unifbd6} 
and \eqref{eq:unifbd11} that 
\begin{equation}\label{eq:Mdebound1}
        \DE{M} \le C +C\int_0^t (\Lebn{\DE{a} }{2} + \lVert \DE{a}
        \rVert_{W^{1,\infty}} + \Lebn{\nabla \DE{v}}{\I})\DE{M}d\tau. 
\end{equation}
We obtain the following Lemma.
\begin{lemma}\label{lem:unifbd}
Let Assumption \ref{asmp:existence} be satisfied with $s>n/2+1$.
There exists $T$ independent of $\delta$ and $\varepsilon$ such that
the solution $(\DE{a},\DE{v})$ is bounded in $C([0,T]; H^s \times
X^{s+1})$ uniformly in $\delta \in ]0,1]$.
\end{lemma}
\begin{proof}
We only estimate the above $\DE{M}(t)$.
Note that $\DE{v}$ vanishes at spatial infinity.
It implies that $\Lebn{\DE{v}}{\I}$ is bounded by $\Sobn{\nabla
  \DE{v}}{s}$ with some constant, since $n\ge 3$.
and $s>n/2+1$. Sobolev embedding and \eqref{eq:Mdebound1}
yield
\begin{equation}\label{eq:Mdebound2}
        \DE{M} (t)\le C +C\int_0^t \(\DE{M}(\tau)\)^{3/2}d\tau.
\end{equation}
Therefore, there exists $\DE{T}>0$ depending only on $\DE{M}(0)$ such
that $\DE{M}(t)$ is bounded by constant times $\DE{M}(0)$ uniformly in
$t \in [0,\DE{T}]$. 
Since $\DE{M}(0)$ is bounded independent of  $\delta$ and $\varepsilon$ by
assumption, $\DE{T}$ can be taken independent of $\delta$ and
$\varepsilon$, as well as the upper bound of $\DE{M}(t)$. 
\end{proof}
\subsection{Existence of the solution to the nonlinear hyperbolic system}
Next we prove the existence of the solution to \eqref{eq:system1}. 
From Lemma~\ref{lem:unifbd}, we see that the sequences
$\{\DE{a}\}_\delta$ and $\{\DE{v}\}_\delta$ are uniformly bounded 
in $C([0,T];H^s)$ and $C([0,T];X^{s+1}\cap L^{\frac{2n}{n-2}})$, respectively.
Therefore, from Ascoli--Arzela's theorem, for a subsequence
$\delta^\prime$ of $\delta$, 
\begin{align*}
        a^\varepsilon_{\delta^\prime} \rightharpoonup a^\varepsilon
        &\wIN C([0,T];H^{s}), \\ 
        v^\varepsilon_{\delta^\prime} \rightharpoonup v^\varepsilon
        &\wIN C([0,T];L^{\frac{2n}{n-2}}), \\ 
        \nabla v^\varepsilon_{\delta^\prime} \rightharpoonup \nabla
        v^\varepsilon &\wIN C([0,T];H^{s}) 
\end{align*}
as $\delta^\prime \to 0$.
Moreover, we have $(a^\varepsilon , v^\varepsilon) \in C_w ([0,T];H^s
\times X^{s+1} \cap L^{\frac{2n}{n-2}})$. 

We shall show that $(a^\varepsilon , v^\varepsilon)$ satisfies
\eqref{eq:system1} 
in $\mathcal{D}^\prime (]0,T]\times \R^n)$.
We fix some $t$. We choose some $s^\prime$ so that $s>s^\prime>n/2+1$.
Then, the above convergences imply
$(a^\varepsilon_{\delta^\prime},\nabla v^\varepsilon_{\delta^\prime})
\to (a^\varepsilon,\nabla v^\varepsilon)$ 
strongly in $C([0,T];H^{s^\prime}_{\text{loc}}\times
H^{s^\prime}_{\text{loc}})$ as $\delta^\prime \to 0$ . 
Then, we deduce that $\DEp{a}$, $\nabla \DEp{a}$, $\DEp{v}$, and
$\nabla \DEp{v}$ converge  uniformly  
in any compact subset of $\R^n$, since $s>n/2+1$ and
$\Lebn{\DEp{v}-v^\varepsilon}{\I}$ is bounded by $\Sobn{\nabla
  \DEp{v}-\nabla v^\varepsilon}{s^\prime}$ 
with some constant. 
Then, we can pass to the limit in all the terms in \eqref{eq:system1},
except possibly the Hartree term. 
Since $(f\ast  g)\ast  h=f\ast  (g\ast  h)$ and $\Jbr{f\ast  g,h} =
\Jbr{f,\check{g}\ast  h}$ with 
$\check{g}(x) = \bar{g}(-x)$, the Hartree term can be rewritten as 
\begin{align*}
        \Jbr{ \lambda \nabla J_{\delta^\prime} (|x|^{-\gamma} \ast |
        \DEp{a}|^2) ,\varphi} 
        = -\lambda \Jbr{ J_{\delta^\prime} | \DEp{a}|^2 ,
        |x|^{-\gamma}\ast  \nabla\varphi }. 
\end{align*}
The function $|x|^{-\gamma}\ast  \nabla\varphi$ is not compactly
supported, but an $\varepsilon/3$-argument shows that the
right hand side tends to 
$-\lambda \Jbr{ |a^\varepsilon|^2 , |x|^{-\gamma}\ast  \nabla\varphi }$.
%
%
Thus, we obtain the solution $(a^\varepsilon,v^\varepsilon)\in C_{w}
([0,T];H^s \times X^{s+1}\cap L^{\frac{2n}{n-2}})$. 
We now claim that this solution is strongly continuous in time.
To prove this, we only have to show that the solution is norm
continuous, that is,  
the function $M^{\varepsilon}(t):=\Sobn{a^\varepsilon(t)}{s}^2 +
\Sobn{\nabla v^\varepsilon(t)}{s}^2$ is continuous in time. 
In the same way as \eqref{eq:Mdebound1}, we have
\begin{equation}\label{eq:completion1}
        \frac{d}{dt} M^{\varepsilon} \le C \(\Lebn{a^\varepsilon}{2} +
        \lVert a^\eps\rVert_{W^{1,\infty}} + \Lebn{\nabla
        v^\varepsilon}{\I}\)M^{\varepsilon}. 
\end{equation}
Since the right hand side is bounded, $M^{\varepsilon}$ is upper
semi-continuous. 
Weak continuities of $a^\varepsilon$ and $v^\varepsilon$ imply the
lower semi-continuity of $M^\varepsilon$. 
Hence, $M^\varepsilon$ is continuous.
\begin{lemma}\label{lem:wexistence}
Let Assumption \ref{asmp:existence} be satisfied.
Suppose $s>n/2+1$.
Let $T$ be given in Lemma~\ref{lem:unifbd}.
For all $\varepsilon \in [0,1]$, 
there exists $(a^\varepsilon,v^\varepsilon)\in C ([0,T];H^s \times
X^{s+1}\cap L^{\frac{2n}{n-2}})$ 
which solves \eqref{eq:system1} in $\mathcal{D}^\prime$.
\end{lemma}
\subsection{Uniqueness}
We next prove the uniqueness of the solution
$(a^\varepsilon,v^\varepsilon)$ by showing 
that if $(a_1^\varepsilon,v^\varepsilon_1)$ and
$(a^\varepsilon_2,v^\varepsilon_2)$ are solutions to
\eqref{eq:system1}, 
in the class $C([0,T];H^s \times X^{s+1}\cap L^{2n/(n-2)})$ for some
$s> n/2+1$, 
then the distance $(a^\varepsilon_1-a^\varepsilon_2, v^\varepsilon_1-
v^\varepsilon_2)$ is equal to zero 
in $L^\I ([0,T];L^2\times \dot{H}^1)$ sense.
Denote $(d^\varepsilon_a,d_v^\varepsilon) :=
(a^\varepsilon_1-a^\varepsilon_2, v^\varepsilon_1-v^\varepsilon_2)$. 
Then, from \eqref{eq:system1}, the system for
$(d_a^\varepsilon,d_v^\varepsilon)$ is rewritten as 
\begin{align}
        &\partial_t d_a^\varepsilon + d_v^\varepsilon\cdot \nabla
        a_1^\varepsilon + v^\varepsilon_2 \cdot \nabla d_a^\varepsilon 
        + \dfrac{1}{2} d_a^\varepsilon\cdot\nabla v^\varepsilon_1 +
        \frac{1}{2}a^\varepsilon_2 \cdot \nabla d_v^\varepsilon 
         =  i\dfrac{\varepsilon}{2}\Delta d_a^\varepsilon ,
        \label{eq:unique1}\\ 
        &\partial_t d_v^\varepsilon + d_v^\varepsilon\cdot\nabla
        v_1^\varepsilon + v^\varepsilon_2 \cdot \nabla d_v^\varepsilon 
        + \lambda\nabla (|x|^{-\gamma}*(d_a^\varepsilon\overline{
        a_1^\varepsilon} +a^\varepsilon_2 \overline{d_a^\varepsilon}))
        = 0.  
\end{align}
Now estimate the $L^2$ norm of $d_a^\varepsilon$.
From the equation in \eqref{eq:unique1}, it holds that
\begin{align*}
        \frac{d}{dt} \Lebn{d_a^\varepsilon}{2}^2 
        &= 2 \Re \Jbr{\partial_t d_a^\varepsilon,d_a^\varepsilon} \\
        \le C &\, |\Re \Jbr{d_v^\varepsilon \cdot \nabla a_1^\varepsilon,
        d_a^\varepsilon}|  
        + C|\Re \Jbr{v^\varepsilon_2 \cdot \nabla d_a^\varepsilon,
        d_a^\varepsilon}| + C|\Re \Jbr{d_a^\varepsilon\cdot\nabla
        v^\varepsilon_1 , d_a^\varepsilon}|\\ 
        &   + C|\Re \Jbr{a^\varepsilon_2 \cdot \nabla d_v^\varepsilon,
        d_a^\varepsilon}| + |\Re \Jbr{i\Delta
        d_a^\varepsilon,d_a^\varepsilon}| . 
\end{align*}
Now, H\"older's inequality and integration by parts show that
\begin{align*}
        |\Re \Jbr{a^\varepsilon_2 \cdot \nabla d_v^\varepsilon,
        d_a^\varepsilon}| 
        &\le \Lebn{a^\varepsilon_2}{\I} \Lebn{\nabla
        d_v^\varepsilon}{2} \Lebn{d_a^\varepsilon}{2}, \\ 
        |\Re \Jbr{v^\varepsilon_2 \cdot \nabla d_a^\varepsilon,
        d_a^\varepsilon}| 
        + |\Re \Jbr{d_a^\varepsilon\cdot\nabla v^\varepsilon_1 ,
        d_a^\varepsilon}| 
        &\le (\Lebn{\nabla v^\varepsilon_1}{\I} + \Lebn{\nabla
        v^\varepsilon_2}{\I}) \Lebn{d_a}{2}^2, \\ 
        \Re \Jbr{i\Delta d_a^\varepsilon,d_a^\varepsilon}
        &=0.
\end{align*}
Another use of H\"older's and Sobolev inequalities shows
\begin{align*}
        |\Jbr{d_v^\varepsilon \cdot \nabla a_1^\varepsilon, d_a^\varepsilon}|
        &\le \Lebn{d_v^\varepsilon}{\frac{2n}{n-2}} \Lebn{\nabla
        a_1^\varepsilon}{n} \Lebn{d_a^\varepsilon}{2}\\ 
        &\le C\Sobn{\nabla a_1^\varepsilon}{\frac{n}{2}-1}
        \Lebn{\nabla d_v^\varepsilon}{2} \Lebn{d_a^\varepsilon}{2} 
\end{align*}
Thus, we end up with the estimate
\begin{align}\label{eq:unique2}
        \frac{d}{dt} \Lebn{d_a^\varepsilon}{2}^2 \le
        C\(\Lebn{d_a^\varepsilon}{2}^2 + \Lebn{\nabla
        d_v^\varepsilon}{2}^2\), 
\end{align}
where the constant $C$ depends on $\Sobn{a_1^\varepsilon}{n/2}$,
$\Lebn{a_2^\varepsilon}{\I}$, and $\Lebn{\nabla v_k^\varepsilon}{2}$
($k=1,2$). 
Similarly, for all $1 \le i,j \le n$, we have the estimates for
$\partial_i d^\varepsilon_{v,j}$: 
\begin{align*}
        \abs{\Jbr{((\partial_i d_v^\varepsilon)\cdot\nabla)
        v^\varepsilon_{1,j}, \partial_i d^\varepsilon_{v,j}} } &\le
        C\Lebn{\nabla v^\varepsilon_{1,j}}{\I} \Lebn{\partial_i
        d_v^\varepsilon}{2}^2, \\ 
        \abs{\Jbr{(d_v^\varepsilon\cdot\nabla) \partial_i
        v^\varepsilon_{1,j},\partial_i d^\varepsilon_{v,j}} } &\le
        C\Lebn{\nabla \partial_i v^\varepsilon_{2,j}}{n}
        \Lebn{d_v^\varepsilon}{\frac{2n}{n-2}}\Lebn{\partial_i
        d^\varepsilon_{v,j}}{2}, \\ 
        \abs{\Jbr{((\partial_i v^\varepsilon_2) \cdot \nabla)
        d^\varepsilon_{v,j},\partial_i d^\varepsilon_{v,j}} } &\le
        C\Lebn{\partial_i v_2^\varepsilon}{\I} \Lebn{\nabla
        d^\varepsilon_{v,j}}{2}^2, \\ 
        \abs{\Jbr{(v^\varepsilon_2 \cdot \nabla) \partial_i
        d^\varepsilon_{v,j},\partial_i d^\varepsilon_{v,j}} } &\le
        C\Lebn{\nabla v_1^\varepsilon}{\I} \Lebn{\partial_i
        d^\varepsilon_{v,j}}{2}^2, \\ 
        \abs{\Jbr{\partial_i\partial_j
        (|x|^{-\gamma}*(d_a^\varepsilon\overline{ a_1^\varepsilon})),
        \partial_id^\varepsilon_{v,j}} } & \le C
        (\Lebn{a_1^\varepsilon}{\I} +
        \Lebn{a_1^\varepsilon}{2})\Lebn{d_a^\varepsilon}{2}
\Lebn{\partial_id^\varepsilon_{v,j}}{2},   \\ 
        \abs{\Jbr{\partial_i\partial_j (|x|^{-\gamma}*(a^\varepsilon_2
        \overline{d_a^\varepsilon})),\partial_i d^\varepsilon_{v,j}} }
        & \le C (\Lebn{a^\varepsilon_2}{\I} +
        \Lebn{a^\varepsilon_2}{2})\Lebn{d_a^\varepsilon}{2}
\Lebn{\partial_id^\varepsilon_{v,j}}{2}, 
\end{align*}
where $v^\varepsilon_{1,j}$ and $ d^\varepsilon_{v,j}$ denote the
$j$-th components of $v^\varepsilon_1$ and $d_v^\varepsilon$,
respectively. 
Summing up over $i$ and $j$, we obtain
\begin{align}\label{eq:unique4}
        \frac{d}{dt} \Lebn{\nabla d_v^\varepsilon}{2}^2 \le
        C\(\Lebn{d_a^\varepsilon}{2}^2 + \Lebn{\nabla
        d_v^\varepsilon}{2}^2\). 
\end{align}
Denote $D(t):=\Lebn{d_a^\varepsilon(t)}{2}^2 + \Lebn{\nabla
  d_v^\varepsilon(t)}{2}^2$: since $D(0)=0$, Gronwall lemma shows that
  $D(t)=0$ for all $t\in [0,T]$.  

\begin{lemma}
The solution $(a^\varepsilon,v^\varepsilon) \in C([0,T];H^s \times
X^{s+1}\cap L^{\frac{2n}{n-2}})$ 
to the system \eqref{eq:system1} given in Lemma \ref{lem:wexistence}
is unique. 
\end{lemma}
\subsection{Completion  of the proof of Proposition \ref{prop:existence-sys}}
Now we complete the proof of existence result.
\begin{proof}[Proof of Proposition \ref{prop:existence-sys}]
We have already shown that the system \eqref{eq:system1} has
a unique solution $(a^\varepsilon , v^\varepsilon) \in C([0,T];H^s
\times X^{s+1}\cap L^{\frac{2n}{n-2}})$. 

%
We show that the existence time $T$ is independent of $s$, thanks to
tame estimates.
In the above proof, the existence time $T$ depends on $s$.
However, once we show the existence of the solution in $[0,T_{s_0}]$
for some $ s_0 > n/2+1$, 
then, for any $s_1>n/2+1$, we deduce from \eqref{eq:completion1} that
\begin{align*}
        M^\varepsilon_{s_1}(t) \le& M^\varepsilon_{s_1}(0) 
        \exp\left(Ct \sup_{0\le \tau \le t}\left(
        \Lebn{a^\varepsilon(\tau)}{2} + \lVert
        a^\varepsilon(\tau)\rVert_{W^{1,\infty}} 
        + \Lebn{\nabla v^\varepsilon(\tau)}{\I}\right)\right) \\
        \le & M^\varepsilon_{s_1}(0) \exp(Ct \sup M^\varepsilon_{s_0}),
\end{align*}
and so that $M^\varepsilon_{s_1}(t)<\infty$ holds for $t\in[0,T_{s_0}]$.
It means that the solution $(a^\varepsilon,v^\varepsilon)$ extends to
time $T_{s_0}$ as a $H^{s_1} \times X^{s_1+1}\cap
L^{\frac{2n}{n-2}}$-valued function, 
that is, $T_{s_1} \ge T_{s_0}$.
The same argument also shows $T_{s_0} \ge T_{s_1}$.
Therefore, $T$ does not depend on $s$.
\end{proof}
\subsection{Construction of $\phi^\varepsilon$}
We finally construct $\phi^\varepsilon$ from $v^\varepsilon$ defined
in Proposition~\ref{prop:existence-sys}. 
Since $v^\eps$ is known, in view of \eqref{eq:originalsystem1}, it is
natural to define 
$\phi^\eps$ as 
\begin{equation}\label{eq:phi}
  \phi^\eps(t,x) = \phi_0(x) -\int_0^t\(\frac{1}{2}\lvert
  v^\eps(\tau,x)\rvert^2 + \l \(|x|^{-\g}\ast  \lvert
  a^\eps\rvert^2\)(\tau,x)\)d\tau. 
\end{equation}
By Assumption~\ref{asmp:existence} and Lemma~\ref{lem:Zhidkov1},
$\phi_0\in L^\infty \cap
L^{\frac{nq_0}{n-q_0}}$. Proposition~\ref{prop:existence-sys} shows
that 
\begin{equation*}
  \lvert
  v^\eps\rvert^2\in C\([0,T]; X^{s+1}\cap
L^{\frac{nq_0}{n-q_0}}\). 
\end{equation*}
Lemma~\ref{lem:Hartree} with $k=2$ and $s_1=s_2=s$ shows that 
\begin{equation*}
\nabla^2\(|x|^{-\g}\ast  \lvert
  a^\eps\rvert^2\) \in C\([0,T];H^s\).
\end{equation*}
By the Hardy--Littlewood--Sobolev inequality, for all $n/(\gamma+1) <
q < \I$, it holds that 
\[
        \Lebn{\nabla(|x|^{-\gamma} \ast  |a^\varepsilon|^2)}{q} \le
        C\Lebn{(|x|^{-\gamma-1} \ast  
        |a^\varepsilon|^2)}{q} \le C\Sobn{a^\varepsilon}{s}^2,
\]
and $\nabla\(|x|^{-\g}\ast  \lvert
  a^\eps\rvert^2\) \in C\([0,T];H^{s+1}\)$.
Moreover, the Sobolev inequality shows $|x|^{-\gamma} \ast 
|a^\varepsilon|^2 \in L^\I$. Therefore, $\phi^\eps$ has the
regularity announced in Theorem~\ref{thm:existence}. To conclude, we
  simply
notice the identity
\begin{equation*}
  \d_t \( \nabla \phi^\eps -v^\eps\)=\nabla \d_t \phi^\eps -\d_t
  v^\eps =0,
\end{equation*}
so that $v^\eps=\nabla \phi^\eps$, and \eqref{eq:phi} yields the
second equation in \eqref{eq:originalsystem1}. 
This completes the proof of Theorem~\ref{thm:existence}.

\section{Asymptotic expansion}
\label{sec:expansion}
\begin{proof}[Proof of Theorem \ref{thm:expansion}]
{\bf First order}.
Suppose that Assumption \ref{asmp:expansion} is satisfied with $N \ge 1$.
We already know that the equation \eqref{eq:system1} has
a unique solution $(a^\varepsilon ,v^\varepsilon) \in C([0,T]; H^s
\times X^{s+1}\cap L^{q_0})$ for all $\varepsilon \in [0,1]$. 
Denote $(b_0,w_0):=(a^\varepsilon,v^\varepsilon)_{|\varepsilon=0}$.
We define
\begin{align}\label{eq:bwL1}
        b^\varepsilon &= \frac{a^\varepsilon - b_0}{\varepsilon}, &
        w^\varepsilon &= \frac{v^\varepsilon - w_0}{\varepsilon}. 
\end{align}
Substituting $a^\varepsilon = b_0 + \varepsilon b^\varepsilon$ and
$v^\varepsilon = w_0 + \varepsilon w^\varepsilon$ into 
the system \eqref{eq:system1}, we obtain the system for
$(b^\varepsilon,w^\varepsilon)$: 
\begin{equation}\label{eq:systemL11}
        \left\{
        \begin{aligned}
                \partial_t b^\varepsilon + w^\varepsilon\cdot \nabla
                b_0 + w_0 \cdot \nabla b^\varepsilon + \frac{1}{2} b_0
                \nabla \cdot w^\varepsilon 
                + \frac{1}{2} b^\varepsilon \nabla \cdot w_0& \\
                +\varepsilon w^\varepsilon \cdot \nabla b^\varepsilon
                + \frac{\varepsilon}{2}b^\varepsilon\nabla\cdot
                w^\varepsilon 
                - i\frac{1}{2}\Delta b_0 &=
                i\frac{\varepsilon}{2}\Delta b^\varepsilon, \\ 
                \partial_t w^\varepsilon + w^\varepsilon \cdot \nabla
                w_0 + w_0 \cdot \nabla w^\varepsilon 
                + \lambda \nabla(|x|^{-\gamma}\ast 2\Re(b_0
                \overline{b^\varepsilon}))& \\ 
                + \varepsilon w^\varepsilon \cdot \nabla w^\varepsilon
                + \lambda\varepsilon \nabla(|x|^{-\gamma}\ast
                |b^\varepsilon|^2) &= 0, 
        \end{aligned}
        \right.
\end{equation}
\begin{align}\label{eq:systemL12}
        b^\varepsilon_{\mid t=0} =& a_1 + \sum_{j=1}^{N-1} \varepsilon^j
        a_{j+1} + \varepsilon^{N-1}r^\varepsilon_N, &
        w^\varepsilon_{\mid t=0} &= 0, 
\end{align}
where we have used the fact that $(b_0,w_0)$ is the solution to the
system \eqref{eq:system1} with $\varepsilon=0$ 
and the assumption that the initial data of $a^\varepsilon$ is written as
$a_0^\varepsilon = a_0 + \varepsilon a_1 + \sum_{j=2}^N \varepsilon^j
a_j + \varepsilon^N r^\varepsilon_N$. 
Since we know that $b^\eps \in C([0,T];H^s)$ and $w^\eps \in
C([0,T]:X^{s+1}\cap L^{\frac{nq_0}{n-q_0}})$ for $\eps>0$,
we just need to prove \emph{a priori} estimates which are independent of $\eps$. 
Mimicking the energy estimates \eqref{eq:unifbd6} and
\eqref{eq:unifbd11}, we obtain 
\begin{multline}\label{eq:boundL1}
        \frac{d}{dt}(\Sobn{b^\varepsilon}{s-2}^2 + \Sobn{\nabla
        w^\varepsilon}{s-2}^2) \\  
        \le C + C(1+ \varepsilon (\Lebn{b^\varepsilon}{2} +
        \lVert b^\varepsilon\rVert_{W^{1,\infty}} + \Lebn{w^\varepsilon}{\I})) 
        (\Sobn{b^\varepsilon}{s-2}^2 + \Sobn{\nabla w^\varepsilon}{s-2}^2),
\end{multline}
where the constant $C$ depends on $\Sobn{b_0}{s}$ and $\Sobn{\nabla w_0}{s}$.
Indeed, the quadratic terms of the system can be handled by the same way
since they are exactly the same as those in the system
\eqref{eq:system1} up to the constant $\varepsilon$. 
We estimate linear terms essentially by the same way.
Note that the integration by parts does not work well, and so that we
need the $H^{s-1}$-boundedness of $b_0$ and $\nabla w_0$. 
The term $i\frac{1}{2}\Delta b_0 $ is also new. By the presence of
this term, $b_0$ is required to be bounded in $H^s$. 

Mimicking the proof of Proposition \ref{prop:existence-sys}, we can
show the existence of a unique solution $(b^\eps,w^\eps) \in C([0,T],H^{s-2}
\times X^{s-1} \cap L^{\frac{2n}{n-2}})$ for all $\eps\in[0,1]$. 
Since $b^\varepsilon(0)$ is uniformly bounded in $H^{s-2}$ by assumption,
we see that the $H^{s-2}\times X^{s-1}$-bound of
$(b^\varepsilon,w^\varepsilon)$ is independent of $\varepsilon$. 
It proves $\Sobn{a^\varepsilon - b_0}{s-2} + \Sobn{\nabla
  v^\varepsilon - \nabla w_0}{s-2}= O(\varepsilon)$. 
Moreover, the existence time $T$ is also independent of $\varepsilon$.
Then, we see from \eqref{eq:bwL1} that, for $\varepsilon>0$,
the existence time for $(b^\varepsilon,w^\varepsilon)$ must be equal
to that for $(a^\varepsilon,v^\varepsilon)$. 
Hence, we conclude that the existence time for $(b^\varepsilon,
w^\varepsilon)$ with $\varepsilon=0$ is also the same. 
Thus, putting $(b_1,w_1) =
(b^\varepsilon,w^\varepsilon)_{|\varepsilon=0}$, we obtain the
solution to the system 
\begin{equation}\label{eq:systembwL11}
        \left\{
        \begin{aligned}
                \partial_t b_1 + w_1 \cdot \nabla b_0 + w_0 \cdot
                \nabla b_1 + \frac{1}{2} b_0 \nabla \cdot w_1 
                + \frac{1}{2} b_1 \nabla \cdot w_0 -
                i\frac{1}{2}\Delta b_0 &=0, \\ 
                \partial_t w_1 + w_1 \cdot \nabla w_0 + w_0 \cdot
                \nabla w_1 + \lambda \nabla(|x|^{-\gamma}\ast 2\Re(b_0
                \overline{b_1})) &= 0, 
        \end{aligned}
        \right.
\end{equation}
\begin{align}\label{eq:systembwL2}
        b_{1\mid t=0} =& a_1 , & w_{1\mid t=0} &= 0.
\end{align}
Since $w_{1\mid t=0} \equiv 0 \in L^r$ for all $r$ and $b_1 \in H^{s-2}$,
we see that $(b_1,w_1) \in C([0,T];H^{s-2}\times X^{s-1}\cap L^{q} )$
for all $q \in ]n/(\gamma+1),\I]$. 
By the similar way as the construction of $\phi^\varepsilon$, we can
construct $\phi_1$ so that $\phi_1 \in L^p$ for all $p \in
]n/\gamma,\I]$ and $\nabla \phi_1 = v_1$. 

{\bf Higher order}.
Let $N \ge 2$ and Assumption \ref{asmp:expansion} be satisfied.
Take $m \in [2,N]$ and assume that, for $1 \le k \le m-1$, the system
\begin{equation}\label{eq:systembwLn1}
        \left\{
        \begin{aligned}
                \partial_t b_k + \sum_{i+j=k} w_i \cdot \nabla b_j  +
                \frac{1}{2} \sum_{i+j=k} b_i \nabla \cdot w_j 
                 - i\frac{1}{2}\Delta b_{k-1} &=0, \\
                \partial_t w_k + \sum_{i+j=k} w_i \cdot \nabla w_j +
                \lambda \sum_{i+j=k} \nabla(|x|^{-\gamma}\ast 2\Re(b_i
                \overline{b_j})) &= 0, 
        \end{aligned}
        \right.
\end{equation}
\begin{align}\label{eq:systembwLn2}
        b_{k\mid t=0} =& a_k , & w_{k\mid t=0} &= 0,
\end{align}
has a unique solution $(b_k,w_k) \in C([0,T]; H^{s-2k} \times
X^{s-2k+1}\cap L^q)$, where $q\in]n/(\gamma+1),\I]$. 
Denote
\begin{align*}
        b^\varepsilon &= \frac{a^\varepsilon
        -\sum_{j=0}^{m-1}\varepsilon^j  b_j}{\varepsilon^m}, & 
        w^\varepsilon &= \frac{v^\varepsilon
        -\sum_{j=0}^{m-1}\varepsilon^j  w_j}{\varepsilon^m}. 
\end{align*}
Then, $(b^\varepsilon,w^\varepsilon)$ satisfies the following system:
\begin{equation}\label{eq:systembwLm1}
        \left\{
        \begin{aligned}
                \partial_t b^\varepsilon
                +\sum_{\ell=0}^{m-1}\varepsilon^\ell & \(  
                        w^\varepsilon\cdot\nabla b_\ell +
                w_\ell\cdot\nabla b^\varepsilon  
                        + \frac{1}{2}b^\varepsilon \nabla \cdot w_\ell
                + \frac{1}{2}b_\ell \nabla \cdot w^\varepsilon\right) \\ 
                &+ \sum_{\ell=0}^{m-1}\varepsilon^\ell
                \sum_{i,j<m, i+j=m+l} \left(w_i \cdot \nabla
                b_j  + \frac{1}{2}  b_i \nabla \cdot w_j\right)  \\ 
                &+ \varepsilon^m w^\varepsilon \cdot\nabla
                b^\varepsilon + \frac{\varepsilon^m}{2}b^\varepsilon
                \nabla \cdot w^\varepsilon 
                - i\frac{1}{2}\Delta b_{m-1} =
                i\frac{\varepsilon}{2}\Delta b^\varepsilon, \\ 
                \partial_t w^\varepsilon 
                + \sum_{\ell=0}^{m-1}&\varepsilon^\ell \left(
                        w^\varepsilon\cdot\nabla w_\ell + w_\ell\cdot\nabla
                w^\varepsilon \right) + 
                        \lambda \nabla(|x|^{-\gamma}\ast 2\Re(b_\ell
                \overline{b^\varepsilon})) \\ 
                &+ \sum_{\ell=0}^{m-1}\varepsilon^\ell
                \sum_{i,j<m,i+j=m+l}  \left(w_i \cdot \nabla
                w_j  
                + \lambda \nabla(|x|^{-\gamma}\ast 2\Re(b_i
                \overline{b_j}))\right) \\ 
                &+ \varepsilon^m w^\varepsilon \cdot \nabla
                w^\varepsilon +  \lambda \varepsilon^m
                \nabla(|x|^{-\gamma}\ast |b^\varepsilon|^2) = 0, 
        \end{aligned}
        \right.
\end{equation}
\begin{align}\label{eq:systembwLm2}
        b^\varepsilon_{\mid t=0} =& a_{m} + \sum_{j=1}^{N-m}
        \varepsilon^j a_{j+m} + \varepsilon^{N-m} r^\varepsilon_N , &
        w^\eps_{\mid t=0} &= 0, 
\end{align}
By  induction on $m$, the system
\eqref{eq:systembwLm1}--\eqref{eq:systembwLm2} 
has a unique solution 
\[
(b^\varepsilon,w^\varepsilon) \in C\([0,T];H^{s-2m}
\times X^{s-2m+1}\)\quad (s-2m+1>n/2).
\]
As the first order, we see that the upper bound of
$(b^\varepsilon,w^\varepsilon)$ and the existence time $T$ is
independent of $\varepsilon$. 
Furthermore, $T$ is the same one as for $(a^\varepsilon,v^\varepsilon)$.
Note that we need the $H^{s-2m+2}$ boundedness of $\Delta b_{m-1}$
to solve the system \eqref{eq:systembwLm1}--\eqref{eq:systembwLm2}.
Denote $(b_m,w_m):=(b^\varepsilon,w^\varepsilon)_{\mid \eps=0}$.
Then, $(b_m,w_m)$ satisfies
\eqref{eq:systembwLn1}--\eqref{eq:systembwLn2} with $k=m$. 
 
Since $w_{m\mid t=0} \equiv 0 \in L^r$ for all $r$ and $b_m \in H^{s-2m}$
with $s-2m > n/2 + 2(N-m) +1 \ge n/2 + 1$, 
we see that $w_m \in C([0,T]; X^{s-2m+1}\cap L^{q} )$ for all $q \in
]n/(\gamma+1),\I]$,  
and that there exists $\varphi_m$ which satisfies $\varphi_m \in
C([0,T],L^p)$ for all $p \in ]n/\gamma,\I]$ and $\nabla \varphi_m =
w_m$. 
\end{proof}
\begin{proof}[Proof of Corollary~\ref{cor:bkw}]
  Corollary~\ref{cor:bkw} is a straightforward consequence of
  Theorem~\ref{thm:expansion}, by considering the asymptotic expansion
  of
  \begin{equation*}
    e^{i\varphi_0/\eps +i\varphi_1}\(b_0 +\eps b_1 +\ldots +\eps^N
    b_N\) e^{i\eps\varphi_2 + i\eps^2 \varphi_3+\ldots
    +i\eps^{N-1}\varphi_{N}}
  \end{equation*}
in powers of $\eps$. The first exponential is not modified, but one
considers the asymptotic expansion of the last two terms in this
product. 
\smallbreak

Note the shift in precision, between
Theorem~\ref{thm:expansion} and Corollary~\ref{cor:bkw}: the initial
order of precision $o(\eps^N)$ becomes $o(\eps^{N-1})$ (in $L^2\cap
L^\infty$). This is 
because the phase $\phi^\eps$ is divided by $\eps$ to go back to
$u^\eps$. This phenomenon has several consequences, see
e.g. \cite{CaARMA} for instabilities. 
\end{proof}
\section{Proof of Corollary~\ref{cor:loss}}
\label{sec:loss}

To see that Corollary~\ref{cor:loss} is a consequence of
Corollary~\ref{cor:bkw}, we resume the same approach as in
\cite{CaARMA}. Let $a_0\in \Sch(\R^n)$ be non-trivial (independent of
$h$), and consider
\begin{equation*}
  \psi_0^h(x)=h^{s-n/2}a_0\(\frac{x}{h}\). 
\end{equation*}
Let $\eps= h^{s_c-s}=h^{\g/2-1-s}$: $h$ and $\eps$ go to zero
simultaneously by assumption. Consider the change of unknown function
\begin{equation*}
  \psi^h(t,x)= h^{s-n/2}u^\eps\(\frac{t}{\eps h^2},\frac{x}{h}\). 
\end{equation*}
Then the Cauchy problem for $\psi^h$ is \emph{equivalent} to:
\begin{equation*}
   i\eps\d_t u^\eps +\frac{\eps^2}{2}\Delta u^\eps= \lambda
                (\lvert x\rvert ^{-\gamma}\ast  |u^\eps|^2)u^\eps\quad
                ;\quad u^\eps_{\mid t=0}=a_0.
\end{equation*}
This is \eqref{eq:r3}--\eqref{eq:ci} with $\phi_0=0$, and $a_0^\eps
=a_0$ independent of $\eps$. By construction, the phase $\varphi_0$
provided in Theorem~\ref{thm:expansion} is such that
\begin{equation*}
  \varphi_{0\mid t=0}=0\quad ;\quad \d_t \varphi_{0\mid t=0} = -\l
  |x|^{-\g}\ast |a_0|^2,
\end{equation*}
where we have used the equation determining $\varphi_0$, that is,
\eqref{eq:originalsystem1} with $\eps=0$. Therefore, there exists
$\tau>0$ independent of $\eps$ such that $\varphi_{0\mid t=\tau}$ is
non-trivial on the support of $a_0$: at time $t=\tau$, $u^\eps$ is
exactly $\eps$-oscillatory, from Corollary~\ref{cor:bkw}. Back to
$\psi^h$, this yields Corollary~\ref{cor:loss}, up to replacing $a_0$
with $|\log h|^{-1}a_0$ (this makes no trouble, since the competition
in terms of $h$ is
logarithmic decay \emph{vs.} algebraic decay).

\subsection*{Acknowledgments} The authors express their deep gratitude
to Professor Yoshio Tsutsumi for his helpful advice.

\bibliographystyle{amsplain}
\bibliography{biblio}

\end{document}